\documentstyle[12pt]{article}
\def\Bbb#1{{\hbox{\bf #1}}}
\def\N{\Bbb N}

\def\Z{\Bbb Z}

\def\T{\Bbb T}
\begin{document}

\title{On the Distribution of Sidon Series}
\author{Nakhl\'e H.\ Asmar\thanks{Supported by NSF grant DMS 9102044}\ \ 
and Stephen Montgomery-Smith\thanks{Supported by NSF grant DMS 9001796}
\\ University of Missouri--Columbia \\ Columbia, MO 65211 \\ U.\ S.\ A.
\date{}}
\maketitle
\begin{abstract}
Let $B$ denote an arbitrary Banach space, $G$\ a compact abelian group
with Haar measure $\mu$\ and
dual group $\Gamma$. Let $E$ be a Sidon subset of $\Gamma$ with Sidon
constant $S(E)$. Let
$r_n$ denote the $n$-th Rademacher function on $[0 , 1]$.
We show that there is a constant $c$, depending only on $S(E)$, such
that, for all $\alpha > 0$:
\begin{eqnarray*}
c^{-1}P\left[\left\| \sum_{n=1}^Na_nr_n\right\| \geq
c\alpha \right] & \leq &
 \mu\left[ \left\|
\sum_{n=1}^Na_n\gamma_n\right\|\geq \alpha \right]
\\
& & \leq \ 
c\,P\left[\left\| \sum_{n=1}^Na_nr_n\right\| \geq c^{-1}\alpha
\right]
\end{eqnarray*}
 where $a_1$, $\ldots$ , $a_N$ are arbitrary elements of $B$,
and $\gamma_1$ ,
$\ldots$ , $\gamma_N$ are arbitrary elements of $E$.  We prove a similar
result for Sidon subsets of dual objects of compact groups,
and apply our results to
obtain new lower bounds for the distribution functions of scalar-valued
Sidon series.
We also note that either one of the above inequalities, even in the scalar
case, characterizes Sidon sets.
\end{abstract}

\pagebreak

\section{Introduction}

Suppose that $G$ is
 a compact abelian group with dual group $\Gamma$.
Denote the normalized Haar measure on $G$ by $\mu$.  Let ${\cal
C}$($G$)
be the Banach space of continuous complex-valued functions on $G$.
 If $S$
$\subset \Gamma$ , a function $f \in L^1(G)$ is called $S-spectral$
whenever $\hat{f}$ is supported in $S$, where here and throughout the
paper
$\hat{}$ denotes taking the Fourier transform.  The collection of
$S$-spectral
functions that belong to a class of functions ${\cal W}$ will be denoted
by
${\cal W}_S$.

\newtheorem{defsidon}{Definition}[section]
\begin{defsidon}
A subset $E$ of $\Gamma$ is called a Sidon set if there is a
constant $c > 0$, depending only on $E$, such that
\begin{equation}
\sum_{\gamma \in \Gamma} \left|\hat{f}(\gamma )\right|
\leq
c\left\| f\right\|_{\infty}
\label{eq1}
\end{equation}
for every $f \in {\cal C}_E(G)$.
The smallest constant c such that (\ref{eq1}\ ) holds is denoted by S(E)
and is called the constant of sidonicity of E, or the Sidon constant of
E.
\label{defsidon}
\end{defsidon}
If $E = \{ \gamma_j\} \subset \Gamma$ is a Sidon set and $\{ a_j\}$ is a
sequence in a Banach space $B$, then the formal series $\sum a_j\gamma_j$
will be referred to as a {\em B-valued Sidon series}.  The norm on a
given Banach space $B$ will be denoted by $\| \cdot \|$, or, sometimes,
by
$\| \cdot \|_B$. \\
It is well-know that Sidon series
share many common properties with Rademacher series.  The following
theorem of Pisier illustrates this fact and will serve as a crucial tool
in our proofs.

\newtheorem{th:Pis1}[defsidon]{Theorem}
\begin{th:Pis1}[{\cite[Th\'eor\`eme 2.1]{bib:Pi1}}]
Suppose that $E = \{ \gamma_n \} \subset \Gamma $ is a Sidon set, that
$B$ is a Banach space, and that $a_1, \ldots, a_N \in B$.  There is a
constant $c_1$, depending only on the Sidon constant $S(E)$, such that,
for every $p \in [1, \infty [$, we have
\begin{equation}
c_1^{-1}\left( {\bf E}\left\| \sum_{n=1}^N
a_nr_n \right\|^p\right)^{\frac{1}{p}}
\leq
\left( \int_G \left\| \sum_{n=1}^N
a_n \gamma_n \right\|^p d\mu \right)^{\frac{1}{p}}
\\
\leq
c_1 \left( {\bf E}\left\| \sum_{n=1}^N
a_nr_n \right\|^p\right)^{\frac{1}{p}}.
\label{ineq:Pis}
\end{equation}
\label{th:Pis1}
\end{th:Pis1}

In view of this similarity between Sidon series and Rademacher series, it
is natural to ask how the distribution function of a Sidon series
compares to the distribution function of a Rademacher series.  Our main
result provides an answer to this question.

\newtheorem{th:main}[defsidon]{Theorem}
\begin{th:main}
Suppose that $E = \{ \gamma_n\} \subset \Gamma$ is a Sidon set, and let
$B$ denote an arbitrary Banach space.  There is a constant $c > 0$ that
depends only on the Sidon constant $S(E)$, such that for all $a_1, \ldots
, a_N \in B$, and all $\alpha > 0$, we have
\begin{equation}
c^{-1}P\left[ \left\| \sum_{n=1}^N a_nr_n\right\|
\geq c\alpha \right]
\leq
\mu\left[ \left\| \sum_{n=1}^N a_n\gamma_n\right\|
\geq \alpha \right]
\leq
cP\left[ \left\| \sum_{n=1}^N a_nr_n\right\|
\geq c^{-1}\alpha \right].
\label{ineq:main}
\end{equation}

\label{th:main}
\end{th:main}

Thus the distribution functions of Sidon and Rademacher series are
equivalent.  Our proof of this result combines well-known properties
of Sidon series with Lemma \ref{lem:main} below.
This Lemma provides sufficient conditions for the
equivalence
of distribution functions.  It applies as well in the setting of
noncommutative groups yielding an analogue of Theorem \ref{th:main}.
Using
the estimates of \cite{bib:Mo1}, we obtain sharp lower
bound
estimates on the distribution of scalar-valued Sidon series on compact
abelian groups.

We should note that the topological implication of our main result is
much easier to show, as is done in \cite{bib:Pi1}, that is, the 
measure topology on the spaces of Rademacher series and the Sidon series
are equivalent.

\section{A principle for the equivalence of distribution functions}

All random variables are defined on some probability space $(\Omega ,
{\cal M}, dP)$.  We denote the set of positive integers by $\N$, the set
of integers by $\Z$, and the circle group by $\T$.
All other notation is as in Section 1.  We start with a
couple of preliminary lemmas.

\newtheorem{lem:1}{Lemma}[section]
\begin{lem:1}
Suppose that $f_1 , f_2, \ldots, f_N$ are independent identically
distributed random variables , and let $f$ be a function with the same
distribution as the $f_j$'s such that
$$
P[ |f| \geq \alpha ] \geq \frac{\theta}{N}
$$
where $\alpha$ and $\theta$ are positive numbers.  Then
$$
P[ \sup_{1\leq j\leq N}|f_j| \geq \alpha ]
\geq
\frac{\theta}{1 + \theta}.$$

\label{lem:1}
\end{lem:1}

{\bf Proof.}  We first show that, for $\theta > 0$ and $N \in \N$,
we have
\begin{equation}
 (1 - \frac{\theta }{N})^N \leq \frac{1}{(1 +
\theta )}.
\label{2.2}
\end{equation}
This follows from the inequalities
$$(1 - \frac{\theta }{N})(1 + \frac{\theta }{N}) \leq 1$$
and
$$(1 + \frac{\theta }{N})^N \geq 1 + \theta .$$
Using (\ref{2.2}\ ) and independence, we get:
\begin{eqnarray*}
P[\sup_j |f_j| \geq \alpha ]     & = &   1 - P[\sup_j |f_j| < \alpha ] \\
                                 & = &   1 -\prod_j P[|f_j| < \alpha ] \\
                                & \geq &  1 - (1 -\frac{\theta}{N})^N \\
                                & \geq &  1 - \frac{1}{(1 + \theta )} \\
                                & = &   \frac{\theta}{(1 + \theta )}.
\end{eqnarray*}

\newtheorem{lem:2}[lem:1]{Lemma}
\begin{lem:2}
Suppose that $f_1 , f_2, \ldots, f_N$ are independent identically
distributed random variables , and let $f$ be a function with the same
distribution as the $f_j$'s such that
$$
P[ |f| \geq \alpha ] \leq \frac{\theta}{N}
$$
where $\alpha$ and $\theta$ are positive numbers.  Then, we have
\begin{equation}
P[ \sup_j|f_j| \geq \alpha ] \leq \theta .
\label{ineq:lem2}
\end{equation}
\label{lem:2}
\end{lem:2}

{\bf Proof.}  We start with the inequality
\begin{eqnarray*}
(1 - x)^n \geq 1 - nx &  & ( 0\le x \le 1),
\end{eqnarray*}
which can be easily proved by using induction and the inequality
\begin{eqnarray*}
(1 - x)(1 - y) \geq 1 - x - y &  & (x , y \geq 0).
\end{eqnarray*}
Hence, for all $0\le \theta \le N$, $N\in\N$, we have

\begin{eqnarray}
(1 - \frac{\theta}{N})^N \geq 1 - \theta &  & (0\le \theta \le N, N\in\N
).
\label{ineq:prf.lem2.1}
\end{eqnarray}
The proof of (\ref{ineq:lem2}) is now straightforward, using
(\ref{ineq:prf.lem2.1}) and independence:

\begin{eqnarray*}
P[\sup_j |f_j| \geq \alpha ]     & = &   1 - P[\sup_j |f_j| < \alpha ] \\
                                 & = &   1 -\prod_j P[|f_j| < \alpha ] \\
                                 & = &  1 - \prod_j(1- P[ |f_j|\geq\alpha
                                                                   ] ) \\
                                & \leq &  1 - (1 - \frac{\theta}{N})^N \\
                                & \leq & \theta .
\end{eqnarray*}

Before stating our main Lemma, we recall two well-known inequalities.
Let $X$ denote a random variable on a probability space $(\Omega, {\cal
M}, dP)$, then, for all $y > 0$, we have:
\begin{eqnarray}
P[|X| \geq y{\bf E}|X|] \leq y^{-1} & & (Chebychev's\ Inequality);
\label{ineq:Cheb}
\end{eqnarray}
and, for $0 < y < 1$, we have
\begin{eqnarray}
P[|X| \geq y{\bf E}|X|] \geq (1 - y)^2\frac{{\bf E}^2|X|}{{\bf E}|X|^2}
                        & & (Paley-Zygmund\ Inequality).
\label{ineq:PZ}
\end{eqnarray}
See \cite[Ineq. II, p. 8]{bib:Ka}.
The statement of our main Lemma now follows. \\

\newtheorem{lem:main}[lem:1]{Lemma}
\begin{lem:main}
Suppose that $X$ and $Y$ are two Banach valued random variables
that are not identically zero, and suppose that $\{
X_n\}$ and $\{ Y_n\}$ are two sequences of independent
Banach valued random variables such that $X_n$ is identically
distributed with $X$, and $Y_n$ is identically distributed with
$Y$.  Suppose that there are constants $c_1$ and $c_2$ such that, for
all positive integers $N$, we have:
\begin{equation}
c_1^{-1} \left\| \sup_{1 \leq j \leq N}\left\| X_j\right\|  \right\|_1
\leq
\left\| \sup_{1 \leq j \leq N}\left\| Y_j\right\|  \right\|_1
\leq
c_1\left\| \sup_{1 \leq j \leq N}\left\| X_j\right\|  \right\|_1 ;
\label{ineq:lem.main.1}
\end{equation}

\begin{equation}
\frac{\left\| \sup_{1 \leq j \leq N}\left\| X_j\right\|  \right\|^2_1}
{\left\| \sup_{1 \leq j \leq N}\left\| X_j\right\|  \right\|^2_2}
\geq c_2 ;
\label{ineq:lem.main.2}
\end{equation}
and
\begin{equation}
\frac{\left\| \sup_{1 \leq j \leq N}\left\| Y_j\right\|  \right\|^2_1}
{\left\| \sup_{1 \leq j \leq N}\left\| Y_j\right\|  \right\|^2_2}
\geq c_2 .
\label{ineq:lem.main.3}
\end{equation}
Then there is a constant $c$, depending only on $c_1$ and $c_2$, such
that, for all $\alpha > 0$, we have:
\begin{equation}
c^{-1}P\left[ \left\| X\right\| \geq c\alpha \right] \leq
P\left[ \left\| Y\right\| \geq \alpha \right] \leq
cP\left[ \left\| X\right\| \geq c^{-1}\alpha \right].
\label{ineq:lem.main.4}
\end{equation}
\label{lem:main}
\end{lem:main}

{\bf Proof.}  We start with the second inequality in
(\ref{ineq:lem.main.4}).  Given an arbitrary $\alpha = \alpha_1 > 0$ with
\begin{equation}
0 < P\left[\left\|Y\right\| \geq \alpha_1\right],
\label{eq:1}
\end{equation}
choose $\nu$ to be the smallest positive integer satisfying:

\begin{equation}
\frac{1}{2\nu}\leq P\left[\left\|Y\right\| \geq \alpha_1\right] \leq
\frac{1}{\nu}.
\label{eq:2}
\end{equation}
From Lemma \ref{lem:1}, it follows that

\begin{equation}
\frac{1}{3} \leq P\left[\sup_{1\leq j\leq \nu} \left\|Y_j\right\| \geq
\alpha_1\right].
\label{eq:3}
\end{equation}
Chebychev's Inequality (\ref{ineq:Cheb}), and (\ref{eq:3}) imply that

\begin{equation}
\frac{1}{3} \leq \frac{1}{\alpha_1}\left\|\sup_{1\leq j\leq \nu}
\left\|Y_j\right\|\right\|_1 .
\label{eq:4}
\end{equation}
From (\ref{ineq:lem.main.1}) and (\ref{eq:4}), we have

\begin{equation}
\frac{1}{3} \leq \frac{c_1}{\alpha_1}\left\|\sup_{1\leq j\leq \nu}
\left\|X_j\right\|\right\|_1 .
\label{eq:5}
\end{equation}
Hence, for any $\alpha_2 > 0$, (\ref{eq:5}) implies that

\begin{equation}
\frac{\alpha_1}{3 \alpha_2 c_1} \leq
\frac{1}{\alpha_2}\left\|\sup_{1\leq j\leq
\nu}
\left\|X_j\right\|\right\|_1.
\label{eq:6}
\end{equation}
In particular, if $\alpha_2 = \frac{\alpha_1}{6c_1}$, we get from
(\ref{eq:6})

\begin{equation}
2\leq
\frac{1}{\alpha_2}\left\|\sup_{1\leq j\leq
\nu}
\left\|X_j\right\|\right\|_1.
\label{eq:7}
\end{equation}
Now we go back to the Paley--Zygmund Inequality (\ref{ineq:PZ}) and apply
it to
\begin{eqnarray*}
\sup_{1\leq j\leq \nu}\left\| X_j\right\| & {\rm with} &
y=\frac{\alpha_2}{\left\| \sup_{1\leq j\leq \nu}\left\| X_j\right\|
\right\|_1} .
\end{eqnarray*}
Taking into account (\ref{ineq:lem.main.2}) and noticing from
(\ref{eq:7}) that $y\leq \frac{1}{2}$, we get

\begin{eqnarray}
\frac{1}{4}c_2 & \leq & P\left[\sup_{1\leq j\leq \nu} \left\| X_j\right\|
\geq \alpha_2\right] \nonumber \\
& = & P\left[\sup_{1\leq j\leq \nu} \left\| X_j\right\|
\geq d\alpha_1\right]
\label{eq:8}
\end{eqnarray}
where $d = \frac{1}{6c_1}$.  Lemma \ref{lem:2} , (\ref{eq:2}), and
(\ref{eq:8}) imply that

\begin{eqnarray}
P\left[\left\| X\right\| \geq d\alpha_1 \right]
                    & \geq & \frac{c_2}{4\nu} \\
                    & \geq & \frac{c_2}{4}P\left[\left\| Y \right\|
                                            \geq \alpha_1\right].
\label{eq:9}
\end{eqnarray}
Take $c^{-1} = \min (d , \frac{c_2}{4})$, then (\ref{eq:9}) shows that

\begin{equation}
P\left[\left\| X\right\| \geq c^{-1}\alpha_1 \right]
                    \geq
c^{-1}P\left[\left\| Y\right\| \geq \alpha_1 \right].
\label{eq:10}
\end{equation}
Note that (\ref{eq:10}) holds for all $\alpha_1$ for which
(\ref{eq:1}) is true.  For all other values of $\alpha_1$, inequality
(\ref{eq:10}) holds trivially.  Thus (\ref{eq:10}) holds for all
$\alpha_1 > 0$.  Now we repeat the proof with $X$ and $Y$ interchanged.
From (\ref{eq:10}) we get:

\begin{equation}
P\left[\left\| Y\right\| \geq c^{-1}\alpha_1 \right]
                    \geq
c^{-1}P\left[\left\| X\right\| \geq \alpha_1 \right],
\label{eq:11}
\end{equation}
for all $\alpha_1 >0$.  Combining (\ref{eq:10}) and (\ref{eq:11}), we
obtain
\begin{equation}
c^{-1}P\left[\left\| X\right\| \geq \alpha_1 \right]
                    \leq
P\left[\left\| Y\right\| \geq c^{-1}\alpha_1 \right]
                    \leq
cP\left[\left\| X\right\| \geq c^{-2}\alpha_1 \right]
\label{eq:12}
\end{equation}
for all $\alpha_1 > 0$.  Equivalently, we have
\begin{equation}
c^{-1}P\left[\left\| X\right\| \geq c\alpha \right]
                    \leq
P\left[\left\| Y\right\| \geq \alpha \right]
                    \leq
cP\left[\left\| X\right\| \geq c^{-1}\alpha \right]
\label{eq:13}
\end{equation}
for all $\alpha > 0$ , which proves (\ref{ineq:lem.main.4}).

\section{The distribution of Banach valued Sidon series}
To prove Theorem \ref{th:main}, Lemma \ref{lem:main} suggests that we
consider independent copies of the given Sidon series.  The
construction of
independent copies of a given trigonometric polynomial on a group $G$ is
easily done on the product group.  The spectra of the resulting
polynomials are supported in a subset of the product of the character
group.  Our first goal in this section is to study properties of this
set.  To simplify the presentation, we will treat the 
commutative and noncommutative cases separately.  Throughout this
section $G$ will denote a compact abelian group with character group
$\Gamma$ and Haar measure $\mu$.  Similar meanings are attributed to
$G_j , \Gamma_j , \mu_j$, respectively.

\newtheorem{join}{Definition}[section]
\begin{join}
Suppose that $E_j\neq \emptyset$ is a subset of $\Gamma_j$, for $j = 1,
\ldots , n$.  The n-fold join of the sets $E_j$ is a subset of
$\prod_{j=1}^{n}\Gamma_j$, denoted by $\bigvee\limits^n_{j=1}E_j$ , and
defined by:
$$
\bigvee ^n_{j=1}E_j= \{ \gamma = (\gamma_1,\ldots , \gamma_n)
\in \prod_{j=1}^n\Gamma_j : \ all\  but\  one\  \gamma_j \in E_j\ are\
0\}
.
$$
If $E\neq\emptyset\subset\Gamma$, the n-fold join of $E$, denoted by
$\bigvee_{j=1}^nE$, is the set $\bigvee_{j=1}^nE_j$, where $E_j=E$ for
all $j=1, \ldots , n$.

\label{def:join}
\end{join}

Thus a generic element $\gamma$ of $\bigvee_{j=1}^nE_j$ is of the form
$\gamma = (0, 0, \ldots , \gamma_j , 0, \ldots , 0)$ where $\gamma_j \in
E_j$.  When $\gamma$ is evaluated at $x=(x_1, \ldots , x_j, \ldots , x_n)
\in G^n$ we get:
$$\gamma (x) = \gamma_j (x_j).$$

Suppose that $S\subset\Gamma$.  For the sake of our proof of Theorem
\ref{th:main},
it turns out that it is sufficient to study the $n$--fold join of the
set $\{1\}\times S\subset \Z\times\Gamma$.  What is needed is the
following simple result.  A more general result concerning joins of Sidon
sets is presented following the proof of Theorem \ref{th:main}.

\newtheorem{lem:sid1}[join]{Lemma}
\begin{lem:sid1}
Suppose $E\subset\Gamma$ is a Sidon set.  Let $T=\{1\}\times
E\subset\Z\times\Gamma$, and let $n$ be an arbitrary positive integer.
The $n$-fold join $\bigvee_{j=1}^n T$ is a Sidon subset of
$\Z^n\times\Gamma^n$ with Sidon constant
$S\left(\bigvee_{j=1}^n T\right)$ equal to $S(E)$.

\label{lem:sid1}
\end{lem:sid1}

\def\scj{$S\left(\bigvee_{j=1}^n T\right)$}
\def\sidon{$\bigvee_{j=1}^n T$}

{\bf Proof.}  It is clear that $S(T)=S(E)$ and that \scj $\geq S(E)$.
Let
$F$ be a trigonometric polynomial with spectrum supported in \scj .
We can write $F$ as:

$$F =\sum_{j=1}^n f_j(t_j,x_j)$$
where
$$ f_j=\sum_{l=1}^{k_j}a_{jl}e^{it_j}\chi_{jl}(x_j)=
 e^{it_j} \sum_{l=1}^{k_j}a_{jl}\chi_{jl}(x_j). $$
For each $j=1,\ldots,n$, pick $x_j$ so that
$$ \left|\sum_{l=1}^{k_j}a_{jl}\chi_{jl}(x_j)\right|=\|f_j\|_\infty,$$
and then pick $t_j$ so that
 $$e^{it_j} \sum_{l=1}^{k_j}a_{jl}\chi_{jl}(x_j)=\|f_j\|_\infty .$$
We have
$$ f_j(t_j,x_j) = \|f_j\|_\infty ,$$
and so
$$F(t_1, t_2, \ldots, t_n, x_1, x_2, \ldots , x_n)=
\sum_{j=1}^n\|f_j\|_\infty\geq\left(S(E)\right)^{-1}
\sum_{j=1}^n\sum_{l=1}^{k_j}|a_{jl}|.$$

Hence (\ref{eq1}) holds for $F$\ with $c = S(E)$, and the proof is complete.

\

We still need one ingredient for the proof of Theorem \ref{th:main}.
This is the Khintchin--Kahane theorem.

\newtheorem{th:KK}[join]{Theorem}
\begin{th:KK}[{\cite[Th\'eor\`eme K]{bib:Pi1}}]
If $0<p<q<\infty$, there is a constant $K_{p,q}$ such that
\begin{equation}
\left( {\bf E}\left\|\sum_{j=1}^na_jr_j\right\|^q\right)^{\frac{1}{q}}
\leq
K_{p,q}\left( {\bf
E}\left\|\sum_{j=1}^na_jr_j\right\|^p\right)^{\frac{1}{p}}
\label{eq:KK}
\end{equation}
for any sequence $\{a_j\}$ in a Banach space $B$

\label{th:KK}
\end{th:KK}

{\bf Proof of Theorem \ref{th:main}}
Let $f(x)= \sum_{n=1}^Na_n\gamma_n(x)\ (x\in G)$ be a trigonometric
polynomial with
$\gamma_n\in S$, and let $Y=e^{it}\sum_{n=1}^Na_n\gamma_n(x)\ $
where $(t,x)\in
\T\times G$. Clearly, $f$ and $Y$ have the same distribution functions.
We apply Lemma \ref{lem:main} with
\begin{eqnarray*}
 X= \sum_{n=1}^Na_nr_n & {\rm and} & Y= e^{it}\sum_{n=1}^Na_n\gamma_n(x).
\end{eqnarray*}
Given $\nu \in \N$, we construct a sequence of independent random
variables $\{Y_j\}$ on $\T^\nu\times G^\nu$, identically distributed with
$Y$, in the
obvious way:  for $(t,x)=(t_1, t_2,\ldots , t_\nu , x_1, \dots, x_\nu)
\in \T^\nu \times G^\nu$, let

\begin{equation}
Y_j(t,x) =Y(t_j,x_j)=e^{it_j}\sum_{n=1}^Na_n\gamma_n(x_j).
\label{eq:yj}
\end{equation}
Write $\gamma_{nj}\ (n = 1, \ldots , N; j=1, \ldots , \nu)$ for the
character in $\Gamma^\nu$ given by:
$$\gamma_{nj}(x)=\gamma_n(x_j)$$
for all $x \in G^\nu$.  Let $\ell^\nu_\infty (B)$ denote the Banach space
consisting of vectors $a=(a_1,\ldots , a_\nu )$,
where $a_j \in B$, equipped with the norm

$$\|a\|_{\ell_\infty^\nu (B)}=\sup_{1\leq j\leq\nu}\left\|a_j\right\| .$$
Let $a_{nj} \in\ell_\infty^\nu (B)$ be the vector whose components are
all zero except the $j-th$ component is equal to $a_n$:
$a_{nj}=a_n\left(\delta_{ij}\right)^\nu_{i=1}\ , n= 1,\ldots , N\ ;
j=1,\ldots , \nu.$  Let
\begin{eqnarray}
{\cal Y}(t, x) = \sum_{j=1}^\nu\sum_{n=1}^Na_{nj}e^{it_j}\gamma_{nj}(x) &
&
((t, x)\in \T^\nu\times G^\nu).
\label{eq:101}
\end{eqnarray}
We clearly have:
\begin{equation}
{\cal Y}(t, x) = \left(Y_1(t, x), Y_2(t, x),\ldots , Y_\nu(t, x)\right)
\label{eq:102}
\end{equation}
for all $(t, x) \in \T^\nu\times G^\nu$.  
For $j = 1,\ldots,\nu$, let
$$X_j(t)=\sum_{n=1}^Na_nr_{nj}(t).$$
Corresponding to ${\cal Y}$,
construct a 
Rademacher sum
${\cal X}$ with values in $\ell_\infty^\nu (B)$:
\begin{eqnarray}
{\cal X}(t) = \sum_{j=1}^\nu\sum_{n=1}^Na_{nj}r_{nj}(t) &  & (t\in
[0 , 1]),
\label{eq:103}
\end{eqnarray}
where ${{(r_{nj})_{n=1}^N}}_{j=1}^\nu$\ is an enumeration of distinct
Rademacher functions, and $a_{nj}$\ is as above.  
Note that ${\cal Y}$ is a Sidon series with
spectrum supported in the join $\bigvee_{j=1}^\nu \{1\}\times E =
\bigvee_{j=1}^n T$ of
$T= \{1\}\times E$. Pisier's Theorem \ref{th:Pis1} implies that
\begin{equation}
c_1^{-1}\left( {\bf E}\sup_{1\leq j\leq\nu}\left\| X_j\right\|\right)
\leq
\int_{G^\nu}\sup_{1\leq j\leq\nu}\left\| Y_j\right\|
\leq
c_1\left( {\bf E}\sup_{1\leq j\leq\nu}\left\| X_j\right\|\right)
\label{eq:104}
\end{equation}
where $c_1$ depends only on $S\left(\bigvee T\right)$, and hence only on
$S(E)$, by Lemma \ref{lem:sid1}.  We have thus obtained
(\ref{ineq:lem.main.1}). It remains
to prove (\ref{ineq:lem.main.2}) and (\ref{ineq:lem.main.3}).  These are
consequences of (\ref{eq:KK}) and Pisier's Theorem \ref{th:Pis1} applied
to the random variables ${\cal X}$ and ${\cal Y}$ above.  Indeed,
applying
(\ref{eq:KK}) with $p=1$ and $q=2$, we obtain:

\begin{equation}
\left( {\bf E}\left(\sup_{1\leq j\leq\nu}\left\|
X_j\right\|\right)^2\right)^{\frac{1}{2}}
\leq
K_{1,2}{\bf E}\sup_{1\leq j\leq\nu}\left\| X_j\right\|
\label{eq:105}
\end{equation}
which proves (\ref{ineq:lem.main.2}).  To get (\ref{ineq:lem.main.3}), we
apply Pisier's Theorem \ref{th:Pis1}, with $p=2$ to the functions ${\cal
X}$ and ${\cal Y}$, then use (\ref{eq:104}) and (\ref{eq:105}) again.

\newtheorem{rem}[join]{Remark}
\begin{rem}
It is worth noting that either of the inequalities (\ref{ineq:main})
characterizes Sidon sets.  For the first inequality, the assertion
follows
from the definition of Sidon sets.  For the second inequality, this
follows by Pisier's characterization of Sidon sets
\cite{bib:Pi3}.

\end{rem}

In our original proof of Theorem \ref{th:main}, we worked directly with
the function $f=\sum a_j\gamma_j$ whose spectrum is supported in a Sidon
set $E$.  The result that we needed concerned the Sidon constant
of the $n$--fold join of $E$.  We present this result in the next
theorem, because of its interest in its own right.  As far as we know,
the best constant in the theorem below is not known.

\newtheorem{th:join}[join]{Theorem}
\begin{th:join}
Suppose that $E$ is a Sidon subset of $\Gamma$ and $n \geq 1$.  Then
$\bigvee_{j=1}^nE$ is a Sidon subset of $\Gamma^n$ with Sidon constant
$S\left( \bigvee_{j=1}^nE\right) \leq 2\pi S(E) + 1$.  In
particular,
$S\left( \bigvee_{j=1}^nE\right)$ is independent of $n$.

\label{th:join}
\end{th:join}

The proof below was kindly communicated to us by D.
Ullrich and other people after him.  It is an easy consequence of the
following lemma which, as D. Ullrich also remarked, may be well-known to
probabilists.

\newtheorem{Ullrich1}[join]{Lemma}
\begin{Ullrich1}
 For $j = 1, \ldots, n$, let $K_j$ denote a compact topological space.
Suppose that $f_j \in {\cal C}(K_j) , j = 1, \ldots , n$.  Suppose
further
that 0 is in the convex hull of $f_j(K_j)$, for each $j$.  Define $f
\in {\cal C} \left( \prod_{j=1}^nK_j\right)$ by $f(x) =
\sum_{j=1}^nf_j(x_j).$ Then
$$\pi\| f\|_\infty \geq \sum_{j=1}^n\|
f_j\|_\infty.$$
Moreover, the constant $\pi$\ is best possible.

\label{lem:Ullrich}
\end{Ullrich1}

{\bf Proof.}  For all $\theta\in [-\pi,\pi[$, and all $x\in\prod_{j=1}^n
K_j$, we
clearly have:

$$\|f\|_\infty\geq\Re \left( e^{i\theta}f(x)\right)=
\sum_{j=1}^n\Re\left( e^{i\theta}f_j(x_j)\right).$$
Choose $b_j \in K_j$ such that $|f_j(b_j)| = \|f_j\|_\infty$.  Write
$|f_j(b_j)|=e^{i\psi_j}f_j(b_j)$ for $\psi_j\in[-\pi,\pi[$.  Now, for each
$\theta\in[-\pi,\pi[$,
$\Re\left( e^{i\theta}f_j(b_j)\right)$
may, or may not, be nonnegative.  If
$\Re\left( e^{i\theta}f_j(b_j)\right)\geq 0$,
set $a_j(\theta ) = b_j$.  If
$\Re\left( e^{i\theta}f_j(b_j)\right) < 0$, define $a_j(\theta )\in K_j$
to be any element of $K_j$ such that
$\Re\left( e^{i\theta}f_j(a_j(\theta))\right)\geq 0.$
This is possible since the convex hull of $f_j(K_j)$ contains
$0$.  Thus, for all $\theta\in[-\pi,\pi[$, we have:
$$\Re\left( e^{i\theta}f_j(a_j(\theta ))\right)
\geq\max \{0 , \cos (\theta - \psi_j)\|f_j\|_\infty\}.$$
Hence,
$$\|f\|_\infty\geq\sum_{j=1}^n \max \{0 , \cos (\theta -
\psi_j)\|f_j\|_\infty\}.$$
Now integrating both sides of the last inequality with respect to
$\theta\in [-\pi , \pi )$ we obtain
$$2\pi\|f\|_\infty\geq 2\sum_{j=1}^n\|f_j\|_\infty$$
from which the desired inequality follows. \\
To see that we have the best constant in the statement of the Lemma,
we consider the following example.  For each $j = 0$, $1,\ldots,$\ $2n-1$,
we let $K_j = \{0,1\}$, say.  Define $f_j(0) = 0$\ and $f_j(1) = e^{ij\frac
\pi n}$.  Clearly we have $\|f_j\|_\infty = 1$.  To compute
$\| f \|_\infty$, notice that 
$$ \| f \|_\infty = \Re\left(e^{i\theta} \sum_{j=0}^{2n-1} f_j(x_j) \right) $$
for some $0\le \theta < 2\pi$ and some $x_j \in K_j$.  Thus it is clear
that 
$$ \| f \|_\infty = \sum_{j=0}^{n-1} \sin(j\frac\pi n + \theta) ,$$
where $0 \le \theta \le \frac \pi n$.  Thus, as $n \to \infty$, 
$$ \| f \|_\infty \approx n \int_0^1 \sin(\pi x) \,dx
                  =  \frac{2n}\pi
                  =  \frac 1\pi \sum_{j=0}^{2n-1} \| f_j\|_\infty .$$

\newtheorem{Ullrich2}[join]{Corollary}
\begin{Ullrich2}
Suppose $E_j\subset\Gamma_j , S(E_j)\leq\kappa$, for each $j=1, \ldots
, N$, and
$0\notin E_j$.  Then $S\left(\bigvee_{j=1}^N E_j\right) \leq
\pi \kappa$.
\label{cor:Ullrich}
\end{Ullrich2}

{\bf Proof.}  Let $E = \bigvee_{j=1}^N E_j$, and suppose that $f\in
{\cal C}_E(\prod_{j=1}^N G_j)$.  Then $f = \sum_{j=1}^N f_j$ with
$f_j
\in
{\cal C}_{E_j}(G_j).$  The fact that 0 $\notin E_j$ shows that
$f_j$ has mean 0, so Lemma \ref{lem:Ullrich} may be applied:

\begin{eqnarray*}
\pi \|f\|_\infty \geq \sum\|f_j\|_\infty
                       & \geq & \kappa^{-1}\sum_{j}
                           \sum_{\gamma\in E_j}|\hat{f_j}(\gamma)| \\
                       & = & \kappa^{-1}
                           \sum_{\gamma\in E}|\hat{f}(\gamma)|.
\end{eqnarray*}

To prove Theorem \ref{th:join}, apply Corollary \ref{cor:Ullrich},
after removing 0 from $E$ and putting it back, if necessary.

The following is a typical application of Theorem \ref{th:main}.
It amounts to transferring, via Theorem \ref{th:main}, a known result
about scalar valued Rademacher series to Sidon series.  We need a definition.

\newtheorem{definition}[join]{Definition}
\begin{definition}
Suppose that $a =\left( a_n\right)_{n=1}^\infty \in \ell^2$.  Define the
sequence $\left( a_n^*\right)_{n=1}^\infty$ to be the nondecreasing
rearrangement of the terms $|a_n|$.  If $0<t<\infty$, we define
$$K_{1,2}(a,t)=\sum_{n=1}^{\left[ t^2\right] }a_n^*
+ t\left( \sum_{\left[ t^2\right] +1}^\infty
(a_n^*)^2\right)^{\frac{1}{2}}$$
where $\left[ t\right]$ denotes the greatest integer part of $t$.
\end{definition}

The next result follows from Theorem \ref{th:main}, and
\cite{bib:Mo1}.
  The latter is the following result for Rademacher series.

\newtheorem{th}[join]{Theorem}
\begin{th}
Suppose that $E=(\gamma_n)\subset\Gamma$ is a Sidon set.  There is a
constant $c>0$ that depends only on $S(E)$ such that for all
$a=(a_n)_{n=1}^\infty\in\ell^2$ we have

$$\mu\left[\left|\sum_{n=1}^\infty a_n\gamma_n\right|
                \geq cK_{1,2}(a,t)\right]
               \leq ce^{-c^{-1}t^2}$$
                 and
$$\mu\left[\left|\sum_{n=1}^\infty a_n\gamma_n\right|
                \geq c^{-1}K_{1,2}(a,t)\right]
                \geq c^{-1}e^{-ct^2}$$
for all $t>0$.
\end{th}

Thus it is possible to calculate rearrangement invariant norms of scalar
valued Sidon series in the spirit of Rodin and Semyonov \cite{bib:RS}.

\def\tr{{\rm tr}\ }

\

We note that these results generalize to noncommutative
compact groups with no further difficulties.
We follow the notation of \cite{bib:Pi2}, Section 5:  $G$ is a compact
group;
$\Sigma$ is the dual object of $G$; $\mu$ is the normalized Haar
measure on $G$.  Thus $\Sigma$ is the set of equivalence classes of the
irreducible representations of $G$.  For each $\iota\in\Sigma$, we let
$U_\iota$ denote a representing element of the equivalence class.  Thus
$U_\iota (x)$ is, for each $x\in G$, a unitary operator on a fixed finite
dimensional Hilbert space $H_\iota$.  The dimension of $H_\iota$ will be
denoted by
$d_\iota$.  
For further details, we refer the reader to \cite{bib:Pi2}, and
\cite{bib:AHA2}.  

In analogy with Definition \ref{defsidon}, a subset
$S\subset\Sigma$ is called a {\em Sidon set} if there is a constant $c$
such that
$$\sum d_\iota\tr |\hat{f}(\iota)|\leq c\| f\|_\infty$$
for every $f\in{\cal C}_S(G)$, where ${\cal C}_S(G)$ is defined as in the 
abelian case.
\def\eps{\varepsilon_\iota}

Let $I$ be a countable indexing set.  An analogue of the Rademacher
functions is defined as a sequence, $\{\eps\}_{\iota\in I}$, of
independent random variables, each $\eps$ being a random $d_\iota\times
d_\iota$ orthogonal matrix, uniformly distributed on the orthogonal group
${\cal O}(d_\iota )$.  These functions are studied in \cite{bib:Pi2} and
\cite{bib:MP}.  Note that we have the analogue of the
Khintchin--Kahane inequality due to 
Pisier--Marcus,
\cite[Corollary 2.12, p. 91]{bib:MP}.

The analogue of Pisier's Theorem \ref{th:Pis1} can be easily
established in this
setting by repeating the proof in \cite{bib:Pi1} and making use of the
properties of Sidon sets on noncommutative groups.  All of these
properties are
found in \cite[Theorem (37.2)]{bib:AHA2}.  For ease of reference, we
state the result below, and omit the proof.

\newtheorem{th:Pis1.noncom}[join]{Theorem}
\begin{th:Pis1.noncom}
Suppose that $P = (\iota) \subset \Sigma $ is a Sidon set, that
$B$ is a Banach space, and that $M_\iota$ is a $d_\iota\times
d_\iota$ matrix with values in $B$. Then there is
a constant $c_1$, depending only on the Sidon constant $S(P)$, such that,
for every $p \in [1, \infty [$, we have
\begin{eqnarray}
c_1^{-1}\left( {\bf E}\left\| \sum_{\iota\in F}
d_\iota\tr (\eps M_\iota ) \right\|^p\right)^{\frac{1}{p}}
& \leq &
\left( \int_G\left\| \sum_{\iota\in F}
d_\iota\tr (U_\iota M_\iota ) \right\|^p\right)^{\frac{1}{p}} \nonumber
\\
& \leq &
c_1\left( {\bf E}\left\| \sum_{\iota\in F}
d_\iota\tr (\eps M_\iota )\right\|^p\right)^{\frac{1}{p}}
\label{ineq:Pis.noncom}
\end{eqnarray}
for any finite subset $F$ of $P$

\label{th:Pis1.noncom}
\end{th:Pis1.noncom}

The result concerning the join of Sidon sets in duals of compact groups
is a straight analogue of Theorem \ref{th:join}.  We omit even the
statement.
We have thus all the necessary ingredients to prove a noncommutative
version of our main Theorem \ref{th:main}.

\newtheorem{th:main.noncom}[join]{Theorem}
\begin{th:main.noncom}
Let $G$ be a compact group with dual object $\Sigma$.\  Suppose that
$P=\{\iota\}\subset\Sigma$ is a Sidon set.  For $n = 1, \ldots, N$, let
$a^n$ denote a $d_{\iota_n} \times d_{\iota_n}$ matrix with entries in a
Banach space $B$.  There is a constant $c>0$ such that, for all $\alpha
>0$, we have:
\begin{eqnarray}
c^{-1}P\left[\left\|\sum_{n=1}^Nd_{\iota_n}\tr (\varepsilon_{\iota_n}
a^n)\right\|
        \geq c\alpha\right]
        & \leq &
\mu\left[\left\|\sum_{n=1}^Nd_{\iota_n}\tr (U_{\iota_n} a^n) \right\|
        \geq \alpha\right] \nonumber \\
      & \leq &
cP\left[\left\|\sum_{n=1}^Nd_{\iota_n}\tr (\varepsilon_{\iota_n}
a^n)\right\|
        \geq c^{-1}\alpha\right].
\label{eq:main.noncom}
\end{eqnarray}
\label{th:main.noncom}
\end{th:main.noncom}

We close our paper by pointing out to the interested reader that,
under suitable conditions, these
methods also apply to commensurate sets of characters of
Pelczy\'nski \cite{bib:Pe}, as indeed they do also to
topological Sidon sets, also described by Pelczy\'nski {\em loc. sit.}


\begin{thebibliography}{Dillo 83}

\bibitem[HR]{bib:AHA2} E. Hewitt and K. A. Ross, ``Abstract Harmonic
Analysis'', Vol. 2, {\bf 152}, Springer--Verlag, Berlin and New York,
1970.

\bibitem[Ka]{bib:Ka} J.P. Kahane,``Some Random Series of Functions'',
Cambridge
Studies in Advanced Mathematics, {\bf 5}, Cambridge University Press,
Cambridge 1985.


\bibitem[MP]{bib:MP}  M. Marcus and G. Pisier, ``Random Fourier Series
with Applications to Harmonic Analysis'', Annals of Mathematics Studies,
Princeton University Press and University of Tokyo Press, Princeton,
1981.

\bibitem[Mo]{bib:Mo1} S. Montgomery-Smith, {\em The distribution of
Rademacher sums}, Proc. AMS, {\bf 109}, (1990), 517--522.

\bibitem[Pe]{bib:Pe} A.~Pelczy\'nski, {\em Commensurate sequences of
characters}, Proc.\ Amer.\ Math.\ Soc.\ {\bf 104}, (1988), 525--531.

\bibitem[Pi 1]{bib:Pi1} G. Pisier, {\em Les in\'egalit\'es de
Kahane--Khintchin d'apr\`es C. Borell}, S\'eminaire sur la g\'eometrie
des \'espaces de Banach, \'Ecole Polytechnique, Palaiseau, Expos\'e
No.~{\bf VII}, (1977--1978).

\bibitem[Pi 2]{bib:Pi2} G. Pisier, {\em De nouvelles caract\'erisations
des
ensembles de Sidon}, Mathematical Analysis and its Applications, Part B,
Edited by L. Nachbin, Advances in Mathematics Supplementary Studies, {\bf
7B}, (1981), 685--726.

\bibitem[Pi 3]{bib:Pi3} G. Pisier, {\em Ensembles de Sidon et processus
gaussiens}, C.\ R.\ Acad.\ Sc. (Paris), S\'erie A, {\bf 286}, (1978),
671--674.

\bibitem[RS]{bib:RS} V.A.~Rodin and E.M.~Semyonov, {\em Rademacher series 
in symmetric
spaces}, Analyse Math.\ {\bf 1}, (1975), 207--222.


\end{thebibliography}
\end{document}